\documentclass[11pt]{amsart}

\newtheorem{theorem}{Theorem}[section]

\newtheorem{prop}[theorem]{Proposition}

\theoremstyle{definition}
\newtheorem{definition}[theorem]{Definition}

\theoremstyle{remark}

\numberwithin{equation}{section}

\newcommand{\ra}{\rightarrow}

\newcommand{\C}{{\mathbb{C}}}

\newcommand{\Ad}{\text{Ad}}

\newcommand{\dist}{\text{dist}}

\newcommand{\lb}{\langle}
\newcommand{\rb}{\rangle}

\newcommand{\mg}{\mathfrak{g}}
\newcommand{\mh}{\mathfrak{h}}
\newcommand{\mk}{\mathfrak{k}}
\newcommand{\mm}{\mathfrak{m}}
\newcommand{\mmp}{\mathfrak{p}}

\newcommand{\ms}{\mathfrak{s}}

\newcommand{\R}{\mathbb{R}}

\newcommand{\bO}{\mathbb{O}}

\newcommand{\la}{\lambda}

\def\so{{\frak {so}}}
\def\su{{\frak {su}}}
\def\u{{\frak u}}
\def\m{{\frak m}}
\def\s{{\frak s}}
\def\be{\begin{equation}}
\def\ee{\end{equation}}
\def\ba{\begin{array}}
\def\ea{\end{array}}
\def\g{{\frak {g}}}
\def\so{{\frak {so}}}
\def\sp{{\frak {sp}}}
\def\su{{\frak {su}}}
\def\u{{\frak u}}
\def\R{{\Bbb R}}
\def\C{{\Bbb C}}
\def\H{{\Bbb H}}
\def\ov{\overline}

\begin{document}

\newcommand{\spacing}[1]{\renewcommand{\baselinestretch}{#1}\large\normalsize}
\spacing{1.14}

\title[Homogeneous metrics]{Homogeneous metrics with nonnegative curvature}

\author{Lorenz Schwachh\"ofer$^\ast$, Kristopher Tapp}

\address{Fakult\"at f\"ur Mathematik\\Technische Universit\"at Dortmund\\Vogelpothsweg 87\\44221 Dortmund\\Germany}
\email{lschwach@math.uni-dortmund.de}

\address{Department of Mathematics\\Saint Joseph University\\5600 City Avenue
Philadelphia, PA 19131}
\email{ktapp@sju.edu}
\thanks{$^\ast$Supported by the Schwerpunktprogramm Differentialgeometrie of the Deutsche
Forschungsgesellschaft}

\begin{abstract}
Given compact Lie groups $H\subset G$, we study the space of $G$-invariant metrics on $G/H$ with nonnegative sectional curvature.  For an intermediate subgroup $K$ between $H$ and $G$, we derive conditions under which enlarging the Lie algebra of $K$ maintains nonnegative curvature on $G/H$.  Such an enlarging is possible if $(K,H)$ is a symmetric pair, which yields many new examples of nonnegatively curved homogeneous metrics.  We provide other examples of spaces $G/H$ with unexpectedly large families of nonnegatively curved homogeneous metrics.
\end{abstract}
\maketitle
Let $H\subset G$ be compact Lie groups, with Lie algebras $\mh\subset\mg$, and let $g_0$ be a bi-invariant metric on $G$.  The space $G/H$ with the induced normal homogeneous metric, denoted $(G,g_0)/H$, has nonnegative sectional curvature.
Little is known about which other $G$-invariant metrics on $G/H$ have nonnegative sectional curvature, except in certain cases.  In all cases where $G/H$ admits a $G$-invariant metric of positive curvature, the problem has been studied along with the determination of which $G$-invariant metric has the best pinching constant; see~\cite{put},\cite{FMV},\cite{VZ}.  When $H$ is trivial, this problem was solved for $G=SO(3)$ and $U(2)$ in~\cite{SMALL}, and partial results for $G=SO(4)$ were obtained in~\cite{HT}.  Henceforth, we identify $G$-invariant metrics on $G/H$ with $\Ad_H$-invariant inner products on $\mmp$ = the $g_0$-orthogonal complement of $\mh$ in $\mg$.

In Section 1, it is an easy application of Cheeger's method to prove that the solution space is star-shaped.  That is, if $g$ is a $G$-invariant metric on $G/H$ with nonnegative curvature, then the inverse-linear path, $g(t)$, from $g(0)=g_0|_\mmp$ to $g(1)=g$ is through nonnegatively curved $G$-invariant metrics.  Here, a path of inner products on $\mmp$ is called ``inverse-linear'' if the inverses of the associated path of symmetric matrices form a straight line.  This observation reduces our problem to an infinitesimal one: first classify the directions, $g'(0)$, one can move away from the normal homogeneous metric such that the inverse-linear path $g(t)$ appears (up to derivative information at $t=0$) to remain nonnegatively curved.  Then, for each candidate direction, check how far nonnegative curvature is maintained along that path.  In Section 2, we derive curvature variation formulas necessary to implement this strategy, inspired by power series derived by M\"uter for curvature along an inverse-linear path~\cite{Muter}.

In Section 3, we consider an intermediate subgroup $K$ between $H$ and $G$, with subalgebra $\mk$, so we have inclusions $\mh\subset\mk\subset\mg$.  Write $\mmp=\mm\oplus\ms$, where $\mm$ is the orthogonal compliment of $\mh$ in $\mk$ and $\ms$ is the orthogonal compliment of $\mk$ in $\mg$.  The inverse-linear path of $G$-invariant metrics on $G/H$ which gradually enlarges $\mk$ is described as follows for all $A,B\in\mmp$:
\begin{equation}\label{sc}g_t(A,B) = \left(\frac{1}{1-t}\right)\cdot g_0(A^\mm,B^\mm)+g_0(A^\ms,B^\ms),\end{equation}
where superscripts denote $g_0$-orthogonal projections onto the corresponding spaces.

This variation scales the fibers of the Riemannian submersion $(G,g_0)/H\ra(G,g_0)/K$.  For $t<0$, these fibers are shrunk, and $g_t$ has nonnegative curvature because it can be redescribed as a submersion metric obtained by a Cheeger deformation:
$$(G/H,g_t) = ((G/H,g_0)\times(K,(-1/t)\cdot g_0))/K.$$
For $t>0$, these fibers are enlarged, and the situation is more complicated.  We will prove:
\begin{theorem}\label{hom}\hspace{1in}
\begin{enumerate}
\item The metric $g_t$ has nonnegative curvature for small $t>0$ if and only if there exists $C>0$ such that for all $X,Y\in\mmp$, $|[X^\mm,Y^\mm]^\mm|\leq C\cdot|[X,Y]|$.
\item In particular, if $(K,H)$ is a symmetric pair, then $g_t$ has nonnegative curvature for small $t>0$, and in fact for all $t\in(-\infty,1/4].$
\end{enumerate}
\end{theorem}

Part 2 provides a large class of new examples of homogeneous metrics with nonnegative curvature.  Notice $t=1/4$ corresponds to the scaling factor $\frac{1}{1-1/4}=4/3$, which appears elsewhere in the literature as an upper limit for enlarging the totally geodesic fibers of certain Riemannian submersions while maintaining nonnegative curvature, including Hopf fibrations~\cite{VZ},\cite{Vo1},\cite{Vo2}, and fibrations of a compact Lie group by cosets of an abelian group~\cite{GZ}.  Wallach proved in~\cite{Wallach} that if $(K,H)$ and $(G,K)$ are rank 1 symmetric pairs and if the triple $(H,K,G)$ satisfies a certain ``fatness'' property,  then the metric $g_t$ has positive curvature for all $t\in(-\infty,1/4)$, $t\neq 0$.  We re-prove Wallach's theorem in Section 3.

When $H$ is trivial, $g_t$ is a left-invariant metric on $G$ scaled up along $\mk$.  Ziller posed the question of when such a metric $g_t$ is nonnegatively curved~\cite{Wolfgang-survey}.  The following answer was found in~\cite{Lorenz}: the metric $g_t$ has nonnegative curvature for small $t>0$ if and only if the semi-simple part of $\mk$ is an ideal of $\mg$; in particular, when $\mg$ is simple, only abelian subalgebras can be enlarged.

When $Ad_H$ acts irreducibly on $\mmp$, there is only a one-parameter family of $G$-invariant metrics on $G/H$ (coming from scaling), all of which are obviously nonnegatively curved.  If there exists an intermediate subalgebra $\mk$, between $\mh$ and $\mg$, then there exists at least a 2-parameter family of $G$-invariant metrics, and many spaces with exactly 2-parameters arise from such an intermediate subalgebra; such spaces were classified in~\cite{DK}.  Thus, Theorem~\ref{hom} addresses the simplest nontrivial case of our classification problem.

Next, in Chapter 4, we show that more arbitrary metric changes preserve nonnegative curvature, assuming a hypothesis which is similar to (but much stronger than) that of Theorem~\ref{hom}:
\begin{theorem}\label{ex}
If there exists $C>0$ such that for all $X,Y\in\mmp$,
$$|X^\mm\wedge Y^\mm|\leq C\cdot|[X,Y]|,$$
then any left invariant metric on $G$ sufficiently close to $g_0$ which is $Ad_H$-invariant and is a constant multiple of $g_0$ on $\ms$ and $\mh$ (but arbitrary on $\mm$) has nonnegative sectional curvature on all planes contained in $\mmp$; hence, the induced metric on $G/H$ has nonnegative sectional curvature. In particular, this hypothesis is satisfied by the following chains $H\subset K\subset G$:
\begin{enumerate}
\item
$Sp(2) \subset SU(4) \subset SU(5)$
\item $SU(3)\subset SU(4)\cong \text{Spin}(6)\subset \text{Spin}(7)$
\item $G_2\subset\text{Spin}(7)\subset\text{Spin}(p+8)$ for $p\in\{0,1\}$, where the second inclusion is the lift of the standard inclusion $SO(7)\subset SO(p+8)$
\item $Spin'(7) \subset Spin(8) \subset Spin(p+9)$ for $p\in\{0,1,2\}$, where $Spin'(7) \subset Spin(8)$ is the image of the spin representation of $Spin(7)$, and the second is again the lift of $SO(8) \subset SO(p+9)$
\item $SU(2)\subset SO(4)\subset G_2$ (Here, $SU(2) \subset SU(3) \subset G_2$, where $SU(3) \subset G_2$ is the isotropy group of $S^6 = G_2/SU(3)$)
\end{enumerate}
\end{theorem}

For the triples above, one is free to choose the initial direction $g'(0)$ of the variation $g(t)$ to be any $\Ad_H$-invariant self-adjoint endomorphism of $\mm$.  For the first, third and fourth triples, the space of such endomorphisms is $1$-dimensional, while for the second triple it is $2$-dimensional. For the fifth triple, the space is $6$-dimensional, but only $3$-dimensional modulo G-equivariant isometry.  In all examples, there is one additional parameter for scaling $\ms$.

Some other spaces are known to admit large-parameter families of nonnegativelty curved homogeneous metrics (\cite{FMV},\cite{put}), but unlike our new examples, these admit positively curved homogeneous metrics.

The statement about nonnegatively curved planes in $G$ is remarkable on its own, since such a metric cannot have nonnegative sectional curvature on {\em all }Êof $G$, unless $\mh$ is abelian~\cite{Lorenz}. Moreover, when constructing nonnegatively curved metrics with normal homogeneous collars, it is precisely the nonnegative curvature of these planes which is needed~\cite{LT2}.


The authors are pleased to thank Wolfgang Ziller for helpful conversations, and the American Institute of Mathematics for hospitality and funding at a workshop on nonnegative curvature in September, 2007, where portions of this work were discussed.
\section{Inverse-linear paths}

In this section, we prove as a quick application of Cheeger's method:

\begin{prop} \label{ILT} If $g$ is a $G$-invariant metric on $G/H$ with nonnegative curvature, then the inverse-linear path, $g(t)$, from $g(0)=g_0|_\mmp$ to $g(1)=g$ is through nonnegatively curved $G$-invariant metrics.
\end{prop}

The case $H=\{e\}$ is found in~\cite{HT}.  We prove this by showing that any $G$-invariant metric with nonnegative curvature on $G/H$ is connected to the normal homogeneous metric $(G,g_0)/H$ via a canonical path of nonnegatively curved $G$-invariant metrics.  See~\cite{HT} for relevant background on Cheeger's method, which is at the heart of the proof.

\begin{proof}[Proof of Proposition~\ref{ILT}]
Let $h$ be an $Ad_H$-invariant inner product on $\mmp$.  Let $M$ denote $G/H$ with the $G$-invariant metric induced by $h$.  Assume that $M$ has nonnegative curvature.  Consider the following family of nonnegatively curved Riemannian submersion metrics on $M$:
$$M_t = \left( M \times\left(G,\frac 1t\cdot g_0\right)\right)/G.$$
Here, $G$ acts diagonally on $M\times G$ as $g\star(p,a) = (g\star p,ag^{-1})$.
This family extends smoothly at $t=0$ to the original metric $M_0=M$.  Notice that each $M_t$ is $G$-invariant, and is therefore induced by some $Ad_H$-invariant inner product, $h_t$, on $\mmp$.  Let $\{e_i\}$ denote a $g_0$-orthonormal basis of $\mmp$ for which $h$ is diagonalized, with eigenvalues $\{\lambda_i\}$.  Then the metrics $h_t$, considered as symmetric matrices with respect to this basis, evolve as follows:
$$h_t = h(I+t\cdot h)^{-1} = \text{diag}\left\{\frac{\lambda_i}{1+t\lambda_i}\right\}.$$
Notice that $M_t$ converges to a point at $t\ra\infty$, but $t\cdot M_t$ converges to the normal homogeneous space $(G,g_0)/H$.

This shows there exists a path of nonnegatively curved $G$-invariant metrics joining $M$ to $(G,g_0)/H$.  We'd like to see that, up to re-parametrization and re-scaling, this path is exactly the inverse-linear path, $\tilde{h}_s$, from $\tilde{h}_0=(G,g_0)/H$ to $\tilde{h}_1=M$.  The initial direction of this inverse-linear path is $\Psi = (I-h^{-1})$, meaning that, in the basis $\{e_i\}$, we have:
$$\tilde{h}_s = (I-s\Psi)^{-1} = \text{diag}\left\{ \frac{1}{1-s(1-\lambda_i^{-1})}   \right\}.$$
It is straightforward now to check that $s\cdot\tilde{h}_s=h_t$ when $s=\frac{1}{1+t}$.
\end{proof}

\section{Curvature variation formulas}
Proposition~\ref{ILT} suggests an infinitesimal strategy for classifying the $G$-invariant metrics with nonnegative curvature on $G/H$.  The first step is to classify the directions, $\Psi$, in which one can move away from a fixed normal homogeneous metric such that curvature variation formulas predict that nonnegative curvature is maintained along the inverse-linear path in that direction.  In this section, we derive the relevant curvature variation formulas.

A path $g_t$ of $\Ad_H$-invariant inner products on $\mmp$ can be described in terms of $g_0|_\mmp$ as:
$$g_t(A,B) = g_0(\Phi_t A,B)$$
for all $A,B\in\mmp$, where $\Phi_t$ is a family of $g_0$-self-adjoint, $\Ad_H$-invariant, positive-definite endomorphims of $\mmp$.  We henceforth assume the path is inverse-linear, which means that $t\mapsto\Phi_t^{-1}$ is linear, so that:
\begin{equation}\label{hi}\Phi_t=(I-t\cdot\Psi)^{-1}\end{equation}
for some $g_0$-self-adjoint, $\Ad_H$-invariant map $\Psi:\mmp\ra\mmp$.  Notice that $\Psi=\frac{d}{dt}|_{t=0}\Phi_t$.

It is useful to henceforth extend $\Phi_t$ and $\Psi$ to be endomorphisms of all of $\mg$ by defining each $\Phi_t$ to be the identity on $\mh$ and defining $\Psi$ to be zero on $\mh$.  Notice that Equation~\ref{hi} still holds for these extensions.

For $X,Y\in\mmp\cong T_{H}(G/H)$, we let $k(t)$ denote the unnormalized sectional curvature with respect to $g_t$ of the vectors $\Phi_t^{-1}X$ and $\Phi_t^{-1}Y$.  The domain of $k$ is the open interval of $t$'s for which $\Phi_t$ represents a non-degenerate metric, which depends on the eigenvalues of $\Psi$.  Notice that $k(0)=0$ if and only if $[X,Y]=0$.  For such initially-zero curvature planes, we will now exhibit a power series expression for $k(t)$.  It is useful to label the following expressions:
\begin{gather*}
A=[\Psi X,Y]+[X,\Psi Y],\\
D_0=[\Psi X,\Psi Y]-\Psi A.
\end{gather*}

\begin{prop}\label{SC} If $X,Y\in\mmp$ commute, then $k(0)=k'(0)=0$,
$k''(0) = \frac 32 |A^\mh|^2$, and
\begin{eqnarray*}
(1/6)k'''(0) & = & \lb A-(3/2)A^\mh,[\Psi X,\Psi Y]\rb
                        +\lb[\Psi X,X],\Psi[\Psi Y,Y]\rb \\
                  &   & -\lb[X,\Psi Y],\Psi A\rb
                        -\lb[\Psi X,Y],\Psi [\Psi X,Y]\rb,\\
\end{eqnarray*}
and for all $t$ in the domain of $k$,
$$k(t) = t^2\cdot(1/2)k''(0) + t^3\cdot(1/6)k'''(0) -\frac 34 t^4\cdot|D_0^\mmp|^2_{g_t}.$$
\end{prop}

\begin{definition} We refer to $\Psi$ (or to the inverse-linear metric variation it determines) as \emph{infinitesimally nonnegative} if for all $X,Y\in\mmp$, there exists $\epsilon>0$ such that $k(t)\geq 0$ for $t\in[0,\epsilon)$.
\end{definition}

This is clearly true for pairs $X,Y$ which don't commute, so it is equivalent to check the condition for pairs which do commute.  This gives:
\begin{prop}
$\Psi$ is infinitesimally nonnegative if and only if for all $X,Y\in\mmp$ such that $[X,Y]=0$ and $A^\mh=0$, we have that $k'''(0)\geq 0$, and $k'''(0)=0$ implies that $D_0^\mmp=0$.
\end{prop}
By Proposition~\ref{ILT}, one will locate all of the nonnegatively curved $G$-invariant metrics on $G/H$ by searching only along infinitesimally nonnegative paths.  This approach was used in~\cite{HT} (in the case where $H$ is trivial) to restrict the space of possible nonnegatively curved left-invariant metrics on $G$.

Proposition~\ref{SC} is a special case of a power-series for $k(t)$, which we now derive,  which does not assume that $X,Y$ commute.  For this general power series, it is useful to denote:
\begin{eqnarray*}
A & = & [\Psi X,Y]+[X,\Psi Y],\\
B & = & [\Psi X,\Psi Y],\\
C & = & [\Psi X,Y]-[X,\Psi Y],\\
D & = & \Psi^2[X,Y] + B - \Psi A
\end{eqnarray*}
With this notation we have:
\begin{prop}\label{gen} For any $X,Y\in\mmp$ and all $t$ in the domain of $k$,
$$k(t) = \alpha + \beta t + \gamma t^2 + \delta t^3 - \frac 34 t^4\cdot|D^\mmp|^2_{g_t}.$$
where
\begin{eqnarray*}
\alpha & = & |[X,Y]^\mh|^2 + \frac 14|[X,Y]^\mmp|^2 \\
\beta & = & -\frac 34\lb \Psi[X,Y],[X,Y]\rb - \frac 32\lb[X,Y]^\mh,A\rb \\
\gamma & = & -\frac 34|\Psi[X,Y]|^2 + \frac 32\lb\Psi[X,Y],A\rb
             -\frac 32\lb[X,Y]^\mm,B\rb +\frac 34|A^\mh|^2 \\
\delta & = &  -\frac 34\lb \Psi^3 [X,Y],[X,Y]\rb+\frac 32 \lb \Psi^2 [X,Y],A \rb
              -\frac 32 \lb \Psi [X,Y],B\rb \\
       &&     -\frac 34 \lb \Psi A,A\rb-\frac 14 \lb \Psi C, C\rb+\lb \Psi[\Psi X,X], [\Psi Y,Y]\rb+\lb A,B\rb -\frac 32\lb A^\mh,B\rb.
\end{eqnarray*}
\end{prop}
\begin{proof}
By O'Neill's formula, $k(t) = \kappa(t)+A(t)$, where $\kappa(t)$ is the unnormalized sectional curvature of $\Phi_t^{-1}X$ and $\Phi_t^{-1}Y$ in the left-invariant metric on $G$ determined by $\Phi_t$, and $A(t)$ is the O'Neill term. Using the expression $\Phi_t^{-1}=I-t\Psi$, we have:
\begin{eqnarray}\label{eA}
\frac 43 A(t)
  & = & |[\Phi_t^{-1}X,\Phi_t^{-1}Y]^\mh|^2 =|[X-t\Psi X,Y-t\Psi Y]^\mh|^2\\
  & = & |[X,Y]^\mh-tA^\mh + t^2B^\mh|^2\notag\\
  & = & |[X,Y]^\mh|^2 -2t\lb[X,Y]^\mh,A\rb
         +t^2\left(|A^\mh|^2 + 2\lb[X,Y]^\mh,B\rb\right) \notag\\
  &&  -2t^3\lb A^\mh,B\rb + t^4|B^\mh|^2.\notag
\end{eqnarray}

It is proven in~\cite{HT} that $\kappa(t) = \overline{\alpha} + \overline{\beta} t + \overline{\gamma}t^2 + \overline{\delta}t^3 - \frac 34 t^4|D|^2_{g_t}$, where
\begin{eqnarray*}
\overline{\alpha} &=& \frac{1}{4} |[X,Y]|^2\\
\overline{\beta} &=& -\frac{3}{4} \lb \Psi [X,Y],[X,Y]\rb\\
\overline{\gamma} &=& -\frac 34 |\Psi[X,Y]|^2+\frac 32 \lb \Psi [X,Y],A\rb
                      -\frac 32 \lb [X,Y],B\rb \\
\overline{\delta} &=& -\frac 34 \lb \Psi^3 [X,Y],[X,Y]\rb+\frac 32 \lb \Psi^2 [X,Y],A \rb
                      - \frac 32 \lb \Psi[X,Y],B\rb \\
                  && - \frac 34 \lb \Psi A,A\rb-\frac 14 \lb \Psi C,
                      C\rb+\lb \Psi[\Psi X,X], [\Psi Y,Y]\rb+\lb A,B\rb.
\end{eqnarray*}
The above expression for $\overline{\gamma}$ is simpler than the one found in~\cite{HT}; to achieve this simplification, use the Jacobi identity to write $\lb[\Psi X,X],[\Psi Y,Y]\rb = \lb[X,Y],B\rb - \lb[X,\Psi Y],[\Psi X,Y]\rb$.

It is straightforward to combine the above power series for $A(t)$ and $\kappa(t)$.  Notice that the $t^4$-term of $k(t)=\kappa(t)+A(t)$ is $\Gamma(t) = \frac 34 t^4(|B^\mh| - |D|^2_{g_t})$, which simplifies because:
$$|B^\mh|^2 - |D|^2_{g_t} = |D^\mh|^2 - (|D^\mh|_{g_t}^2 + |D^\mmp|_{g_t}^2)
                          = |D^\mh|^2 - (|D^\mh|^2 + |D^\mmp|_{g_t}^2) = -|D^\mmp|^2_{g_t}.$$
\end{proof}


\section{Scaling up an intermediate subalgebra}\label{S3}
In this section, we study and prove Theorem~\ref{hom}, which provides conditions under which enlarging an intermediate subalgebra maintains nonnegative curvature on $G/H$.

Suppose $K$ is an intermediate subgroup between $H$ and $G$, with Lie algebra $\mk$, so we have inclusions $\mh\subset\mk\subset\mg$.  Write $\mmp=\mm\oplus\ms$, where $\mm$ is the orthogonal compliment of $\mh$ in $\mk$ and $\ms$ is the orthogonal compliment of $\mk$ in $\mg$.  Let $\Psi$ denote the projecion onto $\mm$, so that $\Psi(A)=A^\mm$ for all $A\in\mg$.  Notice that $\Psi$ determines the inverse-linear path, $g_t$, of $G$-invariant metrics on $G/H$ described in Equation~\ref{sc}, which gradually enlarges the fibers of the Riemannian submersion $(G,g_0)/H\ra(G,g_0)/K$.

We seek conditions under which $g_t$ has nonnegative curvature for small $t>0$.  When $(K,H)$ is a symmetric pair, it is easy to show that $\Psi$ is infinitesimally nonnegative, which provides evidence for Theorem~\ref{hom}.  To fully prove this proposition, we require a power series expression for $k(t)$.

Let $X,Y\in\mmp=\mm\oplus\ms$, and denote
$$M=[X^\mm,Y^\mm],\,\,\,S=[X^\ms,Y^\ms].$$
With this notation, Proposition~\ref{gen} simplifies to:
\begin{eqnarray*}
k(t) & = & \left(\overline{a}|M^\mh|^2+\overline{b}\lb M^\mh,S^\mh\rb+|S^\mh|^2\right) + \left(a|M^\mm|^2+b\lb M^\mm,S^\mm\rb + c|S^\mm|^2\right) +\frac 14|[X,Y]^\ms|^2 \\
     & = & T_1+T_2+T_3.
\end{eqnarray*}
where,
\begin{gather}\label{abc}
\overline{a}=1-3t+3t^2-t^3,\,\,\,\,\,\,\,\,\overline{b}=2-3t,\\
a = \frac 14 - \frac 34\cdot t + \frac 34\cdot t^2 - \frac 14 t^3, \,\,\,\,\,\,\,
b = \frac 12 - \frac 32\cdot t, \,\,\,\,\,\,\,
c = \frac 14-\frac{3t}{4(1-t)}\notag.
\end{gather}

\begin{proof}[Proof of Theorem~\ref{hom}]
Using Cauchy-Swartz, $T_1\geq 0$ when $t\leq 4/3$ because the discriminant is nonnegative:
$$4\overline{a}-\overline{b}^2 = 3t^2-4t^3\geq 0.$$
If $(K,H)$ is a symmetric pair, then $M^\mm=0$, so $T_2 = c|S^\mm|^2$, which is nonnegative for $t\leq 1/4$.  This proves part (2) of the theorem.

For part (1), first assume there exists $C>0$ such that for all $X,Y\in\mmp$, $|M^\mm|\leq C\cdot|[X,Y]|.$  Notice that if $t<1/2$, then
$$T_1\geq \frac {1}{10}|M^\mh+S^\mh|^2 = \frac {1}{10}|[X,Y]^\mh|^2.$$
This is because
$$T_1 - \frac {1}{10}|M^\mh+S^\mh|^2 = \left(\overline{a}-\frac {1}{10}\right)|M^\mh|^2 + \left(\overline{b}-\frac{2}{10}\right)\lb M^\mh,S^\mh\rb + \left(1-\frac{1}{10}\right)|S^\mh|^2,$$
which is nonnegative because the discriminant is nonnegative:
$$\Delta=4\left(\overline{a}-\frac {1}{10}\right)\left(1-\frac {1}{10}\right)-\left(\overline{b}-\frac{2}{10}\right)^2 = \frac 95 t^2 - \frac{18}{5}t^3\geq 0.$$
For $T_2$ we have:
$$T_2\geq a|M^\mm|^2 - b|M^\mm|\cdot|S^\mm| + c|S^\mm|^2 \geq g(t)\cdot|M^\mm|^2,$$
where $g(t) = \frac{4ac-b^2}{4c} = \frac{t^3(t-1)}{1-4t}$.  Notice $g(t)$ is a negative-valued function with $\lim_{t\ra 0} g(t)=0$.

In summary, for $t<1/2$ we have:
\begin{equation}\label{amos}k(t) = T_1+T_2+T_3 \geq \frac {1}{10}|[X,Y]^\mh|^2+g(t)|M^\mm|^2 + \frac 14|[X,Y]^\ms|^2.\end{equation}

At time $t=0$, $T_2 = \frac 14 |M^\mm+S^\mm|^2$, which indicates that for small $t>0$, $T_2$ can only be negative when $M^\mm$ is close to $-S^\mm$.  To make this precise, define ``dist'' as:
$$\text{dist}(A,B) = \max\left\{|\angle(A,B)|,\left|1-|A|/|B|\right|\right\}.$$
Given $\epsilon>0$, we claim there exists $\delta,K>0$ such that if $\text{dist}(M^\mm,-S^\mm)>\epsilon$, then $\frac{T_2}{|M^\mm+S^\mm|^2}\geq K$ for all $t\in[0,\delta]$.  In particular $T_2\geq 0$ (and therefore $k(t)\geq 0$) for $t\in[0,\delta]$.  To see this, notice that $\frac{T_2}{|M^\mm+S^\mm|^2}$ remains unchanged when $M^\mm$ and $S^\mm$ are both scaled by the same factor, so one can assume that the smaller of their lengths equals $1$.  If the larger of their lengths is $\geq 10$, then it is easy to explicitly find $\delta,K$ as above.  When the larger of their lengths is $\leq 10$, a compactness argument suffices to find $\delta,K$.

So it remains to consider the case where $\dist(M^\mm,-S^\mm)<\epsilon$, with $\epsilon>0$ chosen such that
$$ |[X,Y]^\mm|^2 = |M^\mm+S^\mm|^2\leq\frac{1}{2C^2}|M^\mm|^2.$$

In this case, we have by hypothesis:
\begin{eqnarray*}
|M^\mm|^2 & \leq & C^2\cdot\left(|[X,Y]^\mh|^2 + |[X,Y]^\mm|^2 + |[X,Y]^\ms|^2\right)\\
          & \leq & C^2\cdot\left(|[X,Y]^\mh|^2 + \frac{1}{2C^2}|M^\mm|^2 + |[X,Y]^\ms|^2\right).
\end{eqnarray*}
Solving this shows that $|M^\mm|^2\leq 2C^2\left(|[X,Y]^\mh|^2 + |[X,Y]^\ms|^2\right)$.  Combining this with Equation~\ref{amos} shows that $k(t)$ is nonnegative for all $t$ small enough that $2g(t)C^2<1/10$.

The other direction of part (1) of the Theorem follows from similar arguments.
\end{proof}

Next, we recover an important theorem due to Wallach, from which he construct his well-known non-normal homogeneous metrics of positive curvature~\cite{Wallach}.  Recall that the triplet $(H,K,G)$ determines a ``fat homogeneous bundle'' if $[A,B]\neq 0$ for all non-zero $A\in\mm$ and $B\in\ms$; see~\cite{fat} for a survey of literature on fat bundles.
\begin{prop} (Wallach)
If $(K,H)$ and $(G,K)$ are rank 1 symmetric pairs, and $(H,K,G)$ determines a fat homogeneous bundle, then $g_t$ has positive curvature for all $t\in(-\infty,1/4)$, $t\neq 0$.
\end{prop}

\begin{proof}
For linearly independent $X,Y\in\mmp$, if $k(t)=0$ at some non-zero $t\in(-\infty,1/4)$, then the proof of part (2) of Theorem~\ref{hom} implies that $M=[X^\mm,Y^\mm]=0$ and $S=[X^\ms,Y^\ms]=0$ and $[X,Y]^\ms=0$.
So the rank one hypothesis implies that $X^\mm\parallel Y^\mm$ and $X^\ms\parallel Y^\ms$.  Thus, after a change of basis of $\text{span}\{X,Y\}$, we can assume that $X\in\mm$ and $Y\in\ms$.  But then the fact that $[X,Y]^\ms = [X,Y]=0$ contradicts fatness.
\end{proof}

Under the hypotheses of the above proposition, if $k(0)=0$, it is not hard to see that $k''(0)>0$; that is, all initially zero-curvature planes become positively curved to second order.  Since the even derivatives of $k(t)$ are insensitive to the sign of $\Psi$, it does not matter here whether $t$ increases or decreases from zero; in either case, the $A$-tensor makes all initially zero curvature planes become positively curved to second order.

\section{Further Examples}
In this section, we prove Theorem~\ref{ex}, which gives examples of left invariant metrics with many nonnegatively curved planes and, as a consequence, homogeneous spaces with unexpectedly large families of nonnegatively curved homogeneous metrics.  Consider compact Lie groups $H\subset K\subset G$ with Lie algebras $\mh\subset\mk\subset\mg$, and  decompose $\mg = \mh\oplus\mmp = \mh\oplus\mm\oplus\ms$, as in the previous section.

\begin{prop}\label{icecream}
If there exists $C>0$ such that for all $X,Y\in\mmp$,
$$|X^\mm\wedge Y^\mm|\leq C\cdot|[X,Y]|,$$
then any inverse-linear variation $\Phi_t = (I-t\Psi)^{-1}$ of left-invariant $Ad_H$-invariant metrics on $G$ for which $\Psi|_{\ms}=\Psi|_{\mh}=0$ is through metrics which for sufficiently small $t$ have the property that all planes in $\mmp$ are nonnegatively curved.
\end{prop}

The hypothesis of the proposition is clearly stronger than the condition of Theorem~\ref{hom} under which $\mm$ can only be scaled up preserving nonnegative curvature.  Under this stronger hypothesis, the proposition says that arbitrary small changes can be made the metric on $\mm$, and it not only gives information about the metric on $G/H$, but also on $G$.

This proposition clearly implies the first part Theorem~\ref{ex} since any metric close to the normal homogeneous one can be joined by an inverse linear path in that neighborhood.

\begin{proof}
Let $X,Y\in\mmp$.  As in Chapter 2, we have for the curvature of the left-invariant metric on $G$:
$$\kappa(t) = \overline{\alpha} + \overline{\beta} t + \overline{\gamma}t^2 + \overline{\delta}t^3 - \frac 34 t^4|D|^2_{g_t},$$
where the coefficients $\{\overline{\alpha},\overline{\beta},\overline{\gamma},\overline{\delta}\}$ are defined in terms of the expressions $A,B,C,D$.

Notice that $|A^\mk| \leq \lambda_1\cdot |X^\mm\wedge Y^\mm|$, where $\lambda_1$ is the norm of the linear map $\wedge^2\mm\ra\mk$ defined as $x\wedge y\mapsto [\Psi x,y]+[x,\Psi y]$.

Similarly, $|B|\leq \lambda_2\cdot|X^\mm\wedge Y^\mm|$, where $\lambda_2$ is the norm of the linear map $\wedge^2\mm\ra\mk$ defined as $x\wedge y\mapsto [\Psi x,\Psi y]$.

Next, define $E=-\frac 14 \lb \Psi C, C\rb+\lb \Psi[\Psi X,X], [\Psi Y,Y]\rb$, which equals two of the terms in the definition of $\overline{\delta}$.  We claim that $|E|\leq\lambda_3\cdot|X^\mm\wedge Y^\mm|^2$ for some constant $\lambda_3$.  To see this, first consider the symmetric linear map $\rho:\mm\times\mm\ra\mk$ defined as $\rho(x,y) = \frac 12([\Psi x,y] - [x,\Psi y])$.  Next consider the multi-linear map $\Theta:\wedge^2\mm\times\wedge^2\mm\ra\mk$ which is defined as
$$\Theta(x\wedge y,z\wedge w) := \lb \Psi \rho(x,z),\rho(y,w)\rb - \lb\Psi\rho(x,w),\rho(y,z)\rb.$$
Since $\Theta(X^\mm\wedge Y^\mm,X^\mm\wedge Y^\mm)=E$, we may take $\lambda_3$ to be the norm of $\Theta$.

Since the coefficients $\{\overline{\beta},\overline{\gamma},\overline{\delta}\}$ and the term $D$ are defined in terms of the above-bounded expressions, it is a straightforward to use Cauchy-Schwartz to bound their norms and thereby show that there exists a constant $\lambda'$ such that
$$|\overline{\beta}|,|\overline{\gamma}|,|\overline{\delta}|, |D|
     \leq \lambda'\cdot \left(|[X,Y]|^2 + |[X,Y]|\cdot|X^\mm\wedge Y^\mm|
                              + |X^\mm\wedge Y^\mm|^2\right)\leq \lambda\cdot|[X,Y]|^2,$$
where $\lambda = \lambda'(1+C+C^2)$.
In fact, the above bound for $|D|$ also holds for $|D|_{g_t}$ as long at $t$ is small enough that $g_t$ is bounded in terms of $g_0$.
Thus:
\begin{eqnarray*}
\kappa(t) & = & \overline{\alpha} + \overline{\beta} t + \overline{\gamma}t^2 + \overline{\delta}t^3 - \frac 34 t^4|D|^2_{g_t}\\
& \geq & \frac 14|[X,Y]|^2 -(t+t^2+t^3+t^4)\lambda\cdot|[X,Y]|^2
\end{eqnarray*}
which is clearly nonnegative for sufficiently small $t>0$.
\end{proof}


It only remains to prove that the subgroups chains from Theorem~\ref{ex} satisfy the inequality condition of the above proposition.

\begin{prop}
The following triples satisfy the hypothesis of Proposition~\ref{icecream}.
\begin{enumerate}
\item
$Sp(2) \subset SU(4) \subset SU(5)$,
\item
$SU(3) \subset SU(4) \cong Spin(6) \subset Spin(7)$,
\item
$G_2 \subset Spin(7) \subset Spin(p+8)$ for $p\in\{0,1\}$, where the second inclusion is the lift of the inclusion $SO(7) \subset SO(p+8)$.
\item
$Spin'(7) \subset Spin(8) \subset Spin(p+9)$ for $p\in\{0,1,2\}$, where $Spin'(7) \subset Spin(8)$ is the image of the spin representation of $Spin(7)$, and the second is again the lift of $SO(8) \subset SO(p+9)$.
\end{enumerate}
\end{prop}

\begin{proof}
We denote the groups in all cases as $H \subset K \subset G$. Suppose this hypothesis is {\em not} satisfied. Then there exist sequences $\{X_r\}$ and $\{Y_r\}$ in $\m \oplus \s$ such that $X_r^\mm, Y_r^\mm \in \mm$ is an orthonormal pair, and $\lim[X_r, Y_r] = 0$. Passing to a subsequence, we may assume that $X^\mm := \lim X_r^\mm$ and $Y^\mm:= \lim Y_r^\mm$ exist, and we let
\[
B := [X^\mm, Y^\mm] \in \mk.
\]
Since $K/H$ is a sphere and hence the normal homogeneous metric has positive curvature, it follows that $B \neq 0$. Also, $0 = \lim[X_r, Y_r]^\mk = \lim[X_r^\m, Y_r^\m] + [X_r^\s, Y_r^\s]^\mk$, so that
\[
B = -\lim [X_r^\s, Y_r^\s]^\mk,
\]
so that, in particular, we may assume that $[X_r^\s, Y_r^\s]^\mk \neq 0$ for all $r$.

For the first triple, $K/H = Spin(6)/Spin(5) \cong SO(6)/SO(5)$, so that we may regard $X^\m, Y^\m \in \so(5)^\perp \subset \so(6)$, hence $B = [X^\mm, Y^\mm] \in \so(5) \subset \so(6)$ is a matrix of real rank $2$, so that its centralizer is isomorphic to $\so(2) \oplus \so(4)$.

On the other hand, if we regard $[X_r^\s, Y_r^\s]^{\u(4)} \in \u(4) \subset \su(5)$ as a complex matrix where $X_r^\s, Y_r^\s \in \su(4)^\perp \subset \su(5)$, then one verifies that $[X_r^\s, Y_r^\s]^{\u(4)}$ is conjugate to a unique element of the form $diag(\la_1^r i, \la_2^r i, 0, 0)$ with $\la_1^r \geq \la_2^r$. But $\lim [X_r^\s, Y_r^\s]^{\u(4)} = -B \neq 0 \in \su(4)$ exists, so that this limit is conjugate to an element of the form $diag(\la i, -\la i, 0, 0)$ with $\la > 0$, whose centralizer in $\su(4)$ is isomorphic to $\s(\su(2) \oplus \u(1) \oplus \u(1))$. But the centralizer of $B$ is isomorphic to $\so(2) \oplus \so(4)$ which yields the desired contradiction in this case.

For all of the remaining cases we have $G/K = Spin(m)/Spin(n)$ with the inclusion $K \subset G$ induced by the inclusion $SO(n) \subset SO(m)$ for some $(n,m)$.
It follows that for all $X, Y \in \s = \so(n)^\perp$, $[X, Y]^\mk \in \so(n)$ is a matrix which has rank at most $2(m-n)$. Therefore, since $B = -\lim [X_r^\s, Y_r^\s] \neq 0$, it follows that $0 \neq B \in \so(n)$ is a matrix of such a rank.

For the second triple, the rank of $B \in \so(6)$ equals $2(m - n) = 2$, hence its centralizer is isomorphic to $\so(2) \oplus \so(4) \subset \so(6)$.

On the other hand, for $X^\m, Y^\m \in \su(3)^\perp \subset \su(4)$, it is straightforward to verify that $B = [X^\m , Y^\m] \in \su(4)$ is not regular, hence $B$ is conjugate to an element of the form $\text{diag}(\la_1 i, \la_2 i, \la_3 i, 0) \in \su(4)$ with $\la_1 + \la_2 + \la_3 = 0$. Therefore, the centralizer of $B$ is either $\s(\u(1) \oplus \u(1) \oplus \u(1) \oplus \u(1))$ or $\s(\u(2) \oplus \u(1) \oplus \u(1))$, none of which is isomorphic to $\so(2) \oplus \so(4)$ which is a contradiction and finishes the proof for this example.

For the third triple, we will show that for any orthonormal pair $X^\mm, Y^\mm \in \mm$, the rank of $B = [X^\mm, Y^\mm] \in \so(7)$ equals $6$ which will give the desired contradiction as $2(m - n) = 2(p + 1) \leq 4$. For this, we regard $G_2$ as the automorphism group of the octonions $\bO$ which leaves $1 \in \bO$ and hence its orthogonal complement $Im(\bO))$ invariant, and this representation of $G_2$ on $\R^7 \cong Im(\bO)$ lifts to the inclusions $\g_2 \subset \so(Im(\bO))$ and $G_2 \subset Spin(7)$. Then
\[
\so(Im(\bO)) = \g_2 \oplus \{ ad_q: Im(\bO) \longrightarrow Im(\bO)\}
\]
is an orthogonal decomposition, where $ad_q: Im(\bO) \ra Im(\bO)$ is given by $ad_q(x) := q \cdot x - x \cdot q$ since the second summand is $G_2$-equivariantly isomorphic to $Im(\bO)$.

Thus, it remains to show that for an orthonormal pair $q, q' \in Im(\bO)$, the rank of $[ad_q, ad_{q'}] \in \so(Im(\bO))$ equals $6$. Since $G_2$ acts transitively on orthonormal pairs, we may assume that $q = i$ and $q' = j$. Now it is straightforward to verify that the kernel of $[ad_i, ad_j]: Im(\bO) \ra Im(\bO)$ is spanned by $k \in \H$ and is thus one-dimensional.

A similar argument applies to the last case. The orthogonal complement of $\so(7)' \subset \so(8)$ consists of $\{L_q \mid q \in Im(\bO)\}$, where $L_q: \bO \ra \bO$ denotes left multiplication. Assuming w.l.o.g. that $X^\mm = L_i$ and $Y^\mm = L_j$, it is staightforward to verify that $B = [L_i, L_j] \in \so(8)$ is regular, contradicting that $2(m - n) = 2(p + 1) \leq 6$ by assumption.
\end{proof}
\begin{prop}
The triple $SU(2) \subset SO(4) \subset G_2$ satisfies the hypothesis of Lemma~\ref{icecream}.
\end{prop}

\begin{proof}
We decompose the Lie algebra $\g_2$ according to the symmetric pair decomposition of $G_2/SO(4)$ as
\[
\g_2 = (\sp(1)_3 \oplus \sp(1)_1) \oplus \H^2,
\]
where $\sp(1)_3 \subset \sp(2)$ is the Lie algebra spanned by
\[
E_0 := \left(\ba{cc} 3i & \\ & i \ea \right),
E_+:= \left(\ba{cc} 0 & \sqrt 3\\ -\sqrt 3 & 2j \ea \right),
E_-:= \left(\ba{cc} 0 & \sqrt 3 i\\ \sqrt 3 i & 2k \ea \right)
\]
and acts on $\H^2$ from the left, whereas $\sp(1)_1=\text{Im}(\H)$ acts via scalar multiplication from the right. Indeed, one verifies the bracket relations
\[
{}[E_0, E_\pm] = \pm 2 E_\mp, \mbox{ and  } [E_+, E_-] = 2 E_0.
\]

Since $\sp(1)_1$ is the subalgebra which is contained in $\su(3) \subset \g_2$, it follows that in our case, $\m = \sp(1)_3$ and $\s = \H^2$. Thus, we have to show that there cannot be sequences of vectors of the form
\be \label{ex-g2:sequence}
X_n := E_+ + \rho_n \vec v_n \mbox{ and  } Y_n := E_- + \rho'_n \vec w_n
\ee
with unit vectors $\vec v_n, \vec w_n \in \H^2$ and $\rho_n, \rho'_n \geq 0$ such that $\lim[X_n, Y_n] = 0$. By contradiction, we assume that such a sequence of vectors exists and thus may assume that the unit vectors $\vec v := \lim \vec v_n$ and $\vec w := \lim \vec w_n$ exist. Then we have
\[
\ba{lll} 0 & = & \lim \lb E_0, [X_n, Y_n]\rb = \lim \lb[E_0, X_n], Y_n\rb\\
& = & \lim \lb 2 E_- + \rho_n E_0 \cdot \vec v_n, E_- + \rho'_n \vec w_n\rb\\
& = & 2 ||E_-||^2 + \lim \rho_n \rho'_n \lb E_0 \cdot \vec v_n, \vec w_n\rb.
\ea
\]
{}From this we conclude that
\be \label{condition-3}
\liminf \rho_n \rho'_n > 0 \mbox{ and } \lb E_0 \cdot \vec v, \vec w\rb \leq 0.
\ee
Next, for $q \in \sp(1)_1$, we have
\[
0 = \lim \lb[X_n, Y_n], q\rb = \lim \lb X_n, [Y_n, q]\rb = \lim \rho_n \rho'_n \lb\vec v_n, \vec w_n \cdot q\rb,
\]
and since $\liminf \rho_n \rho'_n > 0$, it follows that
\be \label{condition-2}
\lb\vec v, \vec w \cdot q\rb = 0 \mbox{ for all $q \in \sp(1)_1 = \text{Im}(\H)$}.
\ee
Finally,
\[
0 = \lim [X_n, Y_n]^\s = \lim( \rho'_n E_+ \vec w_n - \rho_n E_- \vec v_n).
\]
By (\ref{condition-3}), we may assume that $\rho'_n > 0$ for all $n$. Moreover, $\lim E_- \vec v_n = E_- \vec v \neq 0$ and $\lim E_+ \vec w_n = E_+ \vec w \neq 0$ since $E_\pm$ are regular matrices, so that
\be \label{condition-1}
 0 = E_+ \vec w - c^2 E_- \vec v, \mbox{ where $c^2 := \lim \frac{\rho_n}{\rho'_n} \in (0, \infty)$}.
 \ee

We shall now finish our contradiction by showing that there cannot exist unit vectors $\vec v, \vec w \in \H^2$ satisfying (\ref{condition-3}), (\ref{condition-2}) and (\ref{condition-1}). Namely, $\vec w = c^2 E_+^{-1} E_- \vec v$ by (\ref{condition-1}) , and using the invariance of these conditions under scalar multiplication from the right, we may assume w.l.o.g. that
\[
\vec v = \left( \ba{c} \la\\ z_1 + z_2 j \ea \right), \mbox{ and } \vec w = c^2 E_+^{-1} E_- \vec v = c^2 \left( \ba{c} - \frac 4 {\sqrt 3} \ov z_1 k + \left(\frac 4 {\sqrt 3} \ov z_2 - \la\right) i\\ z_1 i + z_2 k \ea \right),
\]
where $\la \geq 0, c > 0$ and $z_1, z_2 \in \C$. Next, (\ref{condition-2}) holds if for all $q \in \sp(1)_1$,
\[
\ba{lll}
0 & = & \lb\vec v, \vec w \cdot q\rb\\ \\
&  = & c^2\ Re ((\la (- \frac 4 {\sqrt 3} \ov z_1 k + (\frac 4 {\sqrt 3} \ov z_2 - \la) i) + (\ov z_1 - z_2 j)(z_1 i + z_2 k)) q)\\Ê\\
& = & c^2\ Re(((\la(\frac 4 {\sqrt 3} \ov z_2 - \la) + |z_1|^2 - |z_2|^2) i + 2 \ov z_1 (z_2 - \frac 2 {\sqrt 3} \la) k))q).
\ea
\]
If we substitute $q = i$, $q = j$ and $q = k$, we get therefore the equations
\be \label{equation-1}
\la\left(\frac 4 {\sqrt 3} Re(z_2) - \la\right) + |z_1|^2 - |z_2|^2 = 0, \mbox{ and } \overline{z}_1 \left(z_2 - \frac 2 {\sqrt 3} \la\right) = 0.
\ee
If $z_1 \neq 0$, then by the second equation we have $z_2 = \text{Re}(z_2) =  \frac 2 {\sqrt 3} \la$. Substituting this into the first equation of (\ref{equation-1}) implies that $\frac 13 \la^2 + |z_1|^2 = 0$, which is impossible for $z_1 \neq 0$.

Therefore, we conclude from (\ref{equation-1}) that
\be \label{equation-2}
z_1 = 0, \mbox{ and } |z_2|^2 = \la\left(\frac 4 {\sqrt 3} Re(z_2) - \la\right).
\ee
Finally, we calculate from (\ref{equation-2})
\[
\ba{lll}
\left\langle E_0 \cdot \vec v, \vec w\right\rangle & = & c^2 \left\langle E_0 \left( \ba{c} \la\\ z_2 j \ea \right), \left( \ba{c} (\frac 4 {\sqrt 3} \ov z_2 - \la) i\\ z_2 k \ea \right)\right\rangle\\Ê\\
& = & c^2 \left\langle \left( \ba{c} 3 \la i\\ z_2 k \ea \right), \left( \ba{c} (\frac 4 {\sqrt 3} \ov z_2 - \la) i\\ z_2 k \ea \right)\right\rangle\\Ê\\
& = & c^2 (3 \underbrace{\la(\frac 4 {\sqrt 3} Re(z_2) - \la)}_{=|z_2|^2 \mbox{ by (\ref{equation-2})}} + |z_2|^2)\\Ê\\
& = & 4 c^2 |z_2|^2.
\ea
\]
Since $4 c^2 |z_2|^2 = \langle E_0 \cdot \vec v, \vec w\rangle \leq 0$ by (\ref{condition-3}) and $c > 0$, we conclude that $z_2 = 0$, and thus $\la = 0$ by (\ref{equation-2}), i.e. $\vec v = \vec w = 0$. On the other hand, $\vec v$ and $\vec w$ must be unit vectors which is a contradiction.
\end{proof}

\bibliographystyle{amsplain}

\begin{thebibliography}{9}
\bibitem{SMALL} Brown, Finck, Spencer, Tapp, Wu, \emph{Invariant metrics with nonnegative curvature on compact Lie groups}, Cannadian Math. Bull., to appear.
\bibitem{Cheeger} J. Cheeger, \emph{Some examples of manifolds of
nonnegative curvature}, J. Differential Geom. \textbf{8} (1972),
623--628.
\bibitem{DK} W. Dickinson and M. Kerr, \emph{The geometry of compact homogeneous spaces with two isotropy summands}, Ann. Glob. Anal. Geom., to appear
\bibitem{GZ} K.\ Grove and W.\ Ziller, \emph{Curvature and symmetry of Milnor
spheres}, Ann.\ of Math.\ \textbf{152} (2000), 331--367.
\bibitem{HT} J. Huizenga, K. Tapp, \emph{Invariant metrics with nonnegative curvature on $SO(4)$ and other Lie groups}, Michigan Math. J., to appear.
\bibitem{Lorenz} L. Schwachh\"ofer, \emph{A Remark on left invariant metrics on compact Lie groups}, Archiv der Mathematik,  \textbf{90}, No. 2 (2008), 158--162
\bibitem{LT2} L. Schwachh\"ofer, K. Tapp, \emph{Cohomogeneity-one disk bundles with normal homogeneous collars}, preprint.
\bibitem{Muter} M. M\"uter, \emph{Kr\"ummungserh\"ohende Deformationen mittels Gruppenaktionen}, unpublised dissertation, 1987.
\bibitem{put} T. P\"uttmann, \emph{Optimal pinching constants of odd dimensional homogeneous spaces}, Invent. Math. \textbf{138} (1999), no. 3, 631--684.
\bibitem{FMV} F.M. Valiev, \emph{Precise estimates for the sectional curvatures of homogeneous Riemannian metrics on Wallach spaces}, Sib. Math Zhurn. \textbf{200} (1979), 248-262.
\bibitem{VZ} L. Verdianni, W. Ziller, \emph{Positively curved homogeneous metrics on spheres}, Math. Zeit., to appear
\bibitem{Vo1}
D.E.~Volcprime, \emph{ Sectional curvatures of a diagonal family of
${\rm Sp}(n+1)$-invariant metrics of $(4n+3)$-dimensional spheres },
Siberian Math. J. \textbf{35 } (1994), 1089--1100.
\bibitem{Vo2}
D.E.~Volcprime, \emph{ A family of metrics on the $15$-dimensional
sphere}, Siberian Math. J. \textbf{38 } (1997), 223--234.
\bibitem{Wallach} N. Wallach, \emph{Compact Riemannian manifolds with strictly positive curvature}, Ann. of Math. \textbf{96} (1972), 277-295.
\bibitem{fat} W. Ziller, \emph{Fatness revisited}, lecture notes, University of Pennsylvania, 1999.
\bibitem{Wolfgang-survey} W.Ziller, \emph{Examples of Riemannian manifolds with non-negative sectional curvature}, Metric and Comparison Geometry, Surv. Diff. Geom. 11, ed. K.Grove and J.Cheeger, Intern. Press, 2007
\end{thebibliography}

\end{document}